\documentclass[12pt]{article}
\usepackage[english]{babel}
\usepackage{amssymb, amsmath, amsfonts}
\usepackage{mathrsfs}
\newtheorem{Th}{Theorem}
\newtheorem{Lm}{Lemma}

\begin{document}

\begin{center}
{\bf SUBCRITICAL CATALYTIC BRANCHING RANDOM WALK\\
WITH FINITE OR INFINITE VARIANCE OF OFFSPRING NUMBER}
\end{center}
\vskip0,5cm
\begin{center}
Ekaterina Vl. Bulinskaya\footnote{Lomonosov Moscow State
University.}
\end{center}
\vskip1cm

\begin{abstract}
Subcritical catalytic branching random walk on $d$-dimensional
lattice is studied. New theorems concerning the asymptotic behavior
of distributions of local particles numbers are established. To
prove the results different approaches are used including the
connection between fractional moments of random variables and
fractional derivatives of their Laplace transforms. In the previous
papers on this subject only supercritical and critical regimes were
investigated assuming finiteness of the first moment of offspring
number and finiteness of the variance of offspring number,
respectively. In the present paper for the offspring number in
subcritical regime finiteness of the moment of order $1+\delta$
is required where $\delta$ is some positive number.

\vskip0,5cm {\it Keywords and phrases}: branching random walk,
subcritical regime, finite variance of offspring number, infinite
variance of offspring number, local particles numbers, limit
theorems, fractional moments, fractional derivatives.

\vskip0,5cm 2010 {\it AMS classification}: 60F05.
\end{abstract}

\section{Introduction and main results}

The paper is devoted to investigation of \emph{catalytic} branching
random walk (CBRW) on integer lattice $\mathbb{Z}^{d}$,
$d\in\mathbb{N}$. This modification of branching random walk (BRW)
with a single source of branching was proposed by V.A.~Vatutin,
V.A.~Topchij and E.B.~Yarovaya in \cite{VTY} and it comprises
\emph{symmetric} BRW studied earlier (see, e.g.,
\cite{Albeverio_et_al}).

Recall description of CBRW on $\mathbb{Z}^{d}$. Let at the initial
time $t=0$ there be a single particle on the lattice located at
point ${\bf x}\in\mathbb{Z}^{d}$. If ${\bf x}\neq{\bf 0}$ then the
particle performs the continuous time random walk until the first
hitting of the origin. We assume that the random walk outside the
origin is specified by infinitesimal matrix $A=(a({\bf u},{\bf
v}))_{{\bf u},{\bf v}\in\mathbb{Z}^{d}}$ and is symmetric,
homogeneous, having a finite variance of jumps and irreducible (i.e.
the particle passes from an arbitrary point ${\bf u}\in\mathbb{Z}^{d}$
to any point ${\bf v}\in\mathbb{Z}^{d}$ within finite time with
positive probability). It means that
$$a({\bf u},{\bf v})=a({\bf
v},{\bf u}),\quad a({\bf u},{\bf v})=a({\bf 0},{\bf v}-{\bf
u}):=a({\bf v}-{\bf u}),\quad{\bf u},{\bf v}\in\mathbb{Z}^{d},$$
$$\sum\nolimits_{{\bf v}\in\mathbb{Z}^{d}}{a({\bf
v})}=0\quad\mbox{where}\quad a({\bf 0})<0\quad\mbox{and}\quad a({\bf
v})\geq 0\quad\mbox{for}\quad{\bf v}\neq{\bf
0},\quad\sum\nolimits_{{\bf v}\in\mathbb{Z}^{d}}{\|{\bf
v}\|^{2}a({\bf v})}<\infty$$ (here $\|\cdot\|$ denotes an arbitrary
norm in $\mathbb{R}^{d}$). If ${\bf x}={\bf 0}$ or the particle has
just hit the origin then it spends there random time distributed
according to exponential law with parameter $1$. Afterwards it
either dies with probability $\alpha$, producing just before the
death a random offspring number $\xi$ or leaves the origin with
probability $1-\alpha$, so that the intensity of transition from the
origin to point ${\bf v}\neq{\bf 0}$ equals $-(1-\alpha)a({\bf
v})/a({\bf 0}).$ At the origin the branching of particle is
specified by probability generating function
$${f(s):={\sf
E}{s^{\xi}}=\sum\nolimits_{k=0}^{\infty}{f_{k}s^{k}}},\quad
s\in[0,1].$$ At the birth moment the newborn particles are located
at the origin. They evolve in accordance with the scheme described
above independently of each other as well as of their parents history.

The natural objects of study in CBRW are total and local particles
numbers. Denote by $\mu(t)$ the number of particles existing on the
lattice $\mathbb{Z}^{d}$ at time $t\geq0$. In a similar way we
define local numbers $\mu(t;{\bf y})$ as quantities of particles
located at separate points ${\bf y}\in\mathbb{Z}^{d}$ at time $t$.

It was established in \cite{Y_TVP} that exponential growth (as
$t\to\infty$) of both total and local particles numbers in CBRW on
$\mathbb{Z}^{d}$ holds if and only if ${\sf
E}{\xi}>1+h_{d}\alpha^{-1}(1-\alpha)$. Here $h_{d}$ is probability
of the event that a particle which has left the origin will never come
back. Thus, the value ${\sf E}{\xi}=1+h_{d}\alpha^{-1}(1-\alpha)$ is
critical and, similarly to many types of branching processes (see,
e.g., \cite{Sew_74}), CBRW is classified as \emph{supercritical},
\emph{critical} or \emph{subcritical} if the mean offspring number
${\sf E}{\xi}$ is greater, equal or less than
$1+h_{d}\alpha^{-1}(1-\alpha)$, respectively. In view of the random
walk properties such as recurrence or transience, $h_{1}=h_{2}=0$,
whereas $0<h_{d}<1$ for $d\geq3$. Recall that ${\sf E}{\xi}=f'(1)$
and for classic Galton-Watson branching processes the critical value
of $f'(1)$ is equal to 1.

Critical and subcritical CBRW on $\mathbb{Z}^{d}$ are of special
interest since in these cases there arise diverse kinds of limit (in
time) behavior of the total and local particles numbers depending on
dimension $d$. For instance, for $d=1$ and $d=2$ probability ${\sf
P}_{{\bf x}}(\mu(t)>0)$ of non-extinction of the population tends to
zero as time grows (index ${{\bf x}\in\mathbb{Z}^{d}}$ denotes the
starting point of CBRW), whereas for $d\geq3$ this probability has a
positive limit. Such effect is due to existence, on lattice
$\mathbb{Z}^{d}$ for $d\geq3$, of particles which with
positive probability are ever alive and never hit the source
of branching. The total particles number in critical CBRW on integer
line $\mathbb{Z}$, i.e. for $d=1$, was studied in the fundamental paper
\cite{VTY}. In papers \cite{Y_DM} -- \cite{MZ_Bulinskaya} the
investigation was continued for critical and subcritical CBRW on
$\mathbb{Z}^{d}$ for any $d\in\mathbb{N}$.

The analysis of local particles numbers is much more hard. For
critical CBRW on $\mathbb{Z}^{d}$ the limit distributions of local
particles numbers were studied in the series of papers by
V.A.~Vatutin, V.A.~Topchij, Y.~Hu, E.B.~Yarovaya and the author
(see, e.g., \cite{VT_Siberia} -- \cite{DAN_Bulinskaya}). It is
necessary to note that in all these papers the finiteness of second
moment ${\sf E}{\xi}^{2}$ of offspring number was assumed. In the
present paper the limit distributions of local particles numbers in
\emph{subcritical} CBRW on $\mathbb{Z}^{d}$ are studied for the first
time and the corresponding conditions on the moments of offspring
number $\xi$ are less restrictive than those in the critical case.
Namely, we find the asymptotic behavior of the mean local numbers
$m(t;{\bf x},{\bf y}):={\sf E}_{\bf x}{\mu(t;{\bf y})}$ for fixed
${\bf x},{\bf y}\in\mathbb{Z}^{d}$ and $t\to\infty$. Assuming that
${\sf E}{\xi}^{1+\delta}<\infty$ for some $\delta\in(0,1]$, the
similar problem is solved for non-extinction probability $q(t;{\bf
x},{\bf y}):={\sf P}_{\bf x}(\mu(t;{\bf y})>0)$ of local particles
numbers. Moreover, under the same restriction on the moment ${\sf
E}{\xi}^{1+\delta}$ the conditional limit theorems are proved for
$\mu(t;{\bf y})$ as $t\to\infty$.

To formulate the main results let us introduce additional notation.
Let
$$q(s,t;{\bf x},{\bf y}):=1-{\sf
E}_{\bf x}{s^{\mu(t;{\bf y})}}\quad\mbox{and}\quad J(s;{\bf
y}):=\int\nolimits_{0}^{\infty}{\Phi(q(s,t;{\bf 0},{\bf y}))\,dt}$$
where $\Phi(s):=\alpha(f(1-s)-1+f'(1)s)$ and $s\in[0,1]$, $t\geq0$,
${\bf x},{\bf y}\in\mathbb{Z}^{d}$. Consider \emph{transition
probabilities} $p(t;{\bf x},{\bf y}),$ $t\geq0$, ${\bf x},{\bf
y}\in\mathbb{Z}^{d},$ of the random walk generated by matrix~$A$.
According to \cite{VT_Siberia} and Theorem 2.1.1 in
\cite{Yarovaya_book}, for fixed ${\bf x}$ and ${\bf y}$, as
$t\to\infty$, the following asymptotic relations hold true
\begin{equation}\label{p(t;x,y)sim}
p(t;{\bf x},{\bf y})\sim\frac{\gamma_{d}}{t^{d/2}},\quad p\,'(t;{\bf
0},{\bf 0})\sim-\frac{d\,\gamma_{d}}{2\,t^{d/2+1}},\quad p(t;{\bf
0},{\bf 0})-p(t;{\bf x},{\bf y})\sim
\frac{\widetilde{\gamma}_{d}({\bf y}-{\bf x})}{t^{d/2+1}}
\end{equation}
where $\gamma_{d}:=\left((2\pi)^{d}\left|\det{\phi''_{\theta {\bf
\theta}}({\bf 0})}\right|\right)^{-1/2}$, $\phi({\bf
\theta}):=\sum\nolimits_{{\bf z}\in\mathbb{Z}^{d}}{a({\bf 0},{\bf
z})\cos({\bf z},\theta)}$, $\theta\in[-\pi,\pi]^{d},$
$$\phi''_{\theta \theta}({\bf
0}):=\left(\left.\frac{\partial^{2}\phi(\theta)}{\partial \theta_{i}
\partial \theta_{j}}\right|_{\theta={\bf
0}}\right)_{i,j\in\{1,\ldots,d\}},\quad \widetilde{\gamma}_{d}({\bf
z}):=\frac{1}{2(2\pi)^{d}} \int\nolimits_{\mathbb{R}^{d}}{({\bf
v},{\bf z})^{2}e^{\frac{1}{2}(\phi''_{\theta \theta}({\bf 0}){\bf
v},{\bf v})} \,d{\bf v}},\quad {\bf z}\in\mathbb{Z}^{d},$$ and
$(\cdot,\cdot)$ denotes the inner product in $\mathbb{R}^{d}$. Set
$G_{\lambda}({\bf x},{\bf
y}):=\int\nolimits_{0}^{\infty}{e^{-\lambda t}p(t;{\bf x},{\bf
y})\,d t},$ ${\lambda\geq0}$, ${\bf x},{\bf y}\in\mathbb{Z}^{d}$,
i.e. $G_{\lambda}({\bf x},{\bf y})$ is the Laplace transform of the
transition probability $p(\cdot;{\bf x},{\bf y})$. By virtue
of~\eqref{p(t;x,y)sim} for $d\geq3$ the \emph{Green's function}
$G_{0}({\bf x},{\bf y})$ is finite for all ${\bf x},{\bf
y}\in\mathbb{Z}^{d}$, which means the transience of the random walk.
However, for $d=1$ or $d=2$, we deal with recurrent random walk
since $\lim\nolimits_{\lambda\to0+}{G_{\lambda}({\bf x},{\bf
y})}=\infty$. In view of the same formula \eqref{p(t;x,y)sim}
the function $\lim\nolimits_{\lambda\to0+}{(G_{\lambda}({\bf 0},{\bf
0})-G_{\lambda}({\bf x},{\bf y}))}$ is finite for all
$d\in\mathbb{N}$ and ${\bf x},{\bf y}\in\mathbb{Z}^{d}$. So we may
introduce the function
$$\rho_{d}({\bf
z}):=\left\{
\begin{array}{lcl}
(1-\alpha)a^{-1}-\beta\int\nolimits_{0}^{\infty}{\left(p(t;{\bf
0},{\bf 0})-p(t;{\bf 0},{\bf z})\right)dt}&\mbox{if}&\quad
{\bf z}\neq{\bf 0},\\
1&\mbox{if}&\quad{\bf z}={\bf 0}
\end{array}
\right.
$$
where for the sake of convenience we set $a:=-a({\bf 0})$ and
$\beta:=\alpha(f'(1)-1)$. As shown in \cite{VT_Siberia}, $h_{d}=(a
G_{0}({\bf 0},{\bf 0}))^{-1}$. It implies that in subcritical regime
${\beta<(1-\alpha)(a G_{0}({\bf 0},{\bf 0}))^{-1}}$ and, therefore, $\rho_{d}(\cdot)$ is a strictly positive function for all
$d\in\mathbb{N}$. Note also that in accordance with formula
(2.1.15) in \cite{Yarovaya_book} inequality $p(t;{\bf 0},{\bf
0})\geq p(t;{\bf x},{\bf y})$, $t\geq0$, is valid and the function
$p(\cdot;{\bf x},{\bf y})$ is symmetric and homogeneous in ${\bf
x},{\bf y}\in\mathbb{Z}^{d}$.

For $d=1$ and ${\bf x},{\bf y}\in\mathbb{Z}$, we define
$$C_{1}({\bf 0},{\bf y}):=\frac{1-\alpha}{2a\gamma_{1}\pi\beta^{2}}\,\rho_{1}({\bf
y}),\quad C_{1}({\bf x},{\bf
y}):=\frac{1}{2\gamma_{1}\pi\beta^{2}}\,\rho_{1}({\bf
x})\,\rho_{1}({\bf y})+\widetilde{\gamma}_{1}({\bf
x})+\widetilde{\gamma}_{1}({\bf y})-\widetilde{\gamma}_{1}({\bf
y}-{\bf x}),\quad {\bf x}\neq{\bf 0}.$$ In a similar way, for $d=2$
and ${\bf x},{\bf y}\in\mathbb{Z}^{2}$ set
$$C_{2}({\bf 0},{\bf y}):=\frac{1-\alpha}{a\gamma_{2}\beta^{2}}\rho_{2}({\bf
y}),\quad C_{2}({\bf x},{\bf
y}):=\frac{1}{\gamma_{2}\beta^{2}}\rho_{2}({\bf x})\rho_{2}({\bf
y}),\quad {\bf x}\neq{\bf 0}.$$ At last, for $d\geq3$, specify
functions $C_{d}({\bf x},{\bf y})$, ${\bf x},{\bf
y}\in\mathbb{Z}^{d}$, by way of
$$C_{d}({\bf 0},{\bf y}):=\frac{(1-\alpha)\,a\,\gamma_{d}}{(1-\alpha-a\,\beta\,G_{0}({\bf 0},{\bf 0}))^{2}}\,
\rho_{d}({\bf y}),\;\; C_{d}({\bf x},{\bf
y}):=\frac{a^{2}\,\gamma_{d}}{(1-\alpha-a\,\beta\,G_{0}({\bf 0},{\bf
0}))^{2}}\,\rho_{d}({\bf x})\,\rho_{d}({\bf y}),\;\; {\bf x}\neq{\bf
0}.$$

\begin{Th}\label{T-1}
Let ${\sf E}{\xi}<1+h_{d}\alpha^{-1}(1-\alpha)$. Then, as
$t\to\infty$, for each ${\bf x},{\bf y}\in\mathbb{Z}^{d}$, the following
relations are true
\begin{eqnarray*}
m(t;{\bf x},{\bf y})\sim\frac{C_{1}({\bf x},{\bf y})}{t^{3/2}},\quad &d=1,&\\
m(t;{\bf x},{\bf y})\sim\frac{C_{2}({\bf x},{\bf y})}{t\ln^{2}{t}},\quad &d=2,&\\
m(t;{\bf x},{\bf y})\sim\frac{C_{d}({\bf x},{\bf y})}{t^{d/2}},\quad
&d\geq3,&
\end{eqnarray*}
and the functions $C_{d}(\cdot,\cdot)$, $d\in\mathbb{N}$, introduced
above, are strictly positive.
\end{Th}

The statement of Theorem \ref{T-1} generalizes the corresponding
results of Chapter 5 in \cite{Yarovaya_book} concerning the
asymptotic behavior of the first moments of local particles numbers
in subcritical \emph{symmetric} BRW on $\mathbb{Z}^{d}$. Here, in contrast
to \cite{Yarovaya_book}, we use the approaches such as representation
of complex-valued measures in terms of Banach algebras (see
\cite{VT_Siberia}) and Tauberian theorems for derivatives of Laplace
transform (see Section 7.3 in \cite{V_Lectures}).

\begin{Th}\label{T-2}
If ${\sf E}{\xi}<1+h_{d}\alpha^{-1}(1-\alpha)$ and there exists
${\sf E}{\xi^{1+\delta}}$ for some $\delta\in(0,1]$ then, for fixed
${\bf x},{\bf y}\in\mathbb{Z}^{d}$, as $t\to\infty$, one has
\begin{eqnarray*}
q(t;{\bf x},{\bf y})\sim\frac{C_{1}({\bf x},{\bf y})-C_{1}({\bf x},{\bf 0})J(0;{\bf y})}{t^{3/2}},\quad &d=1,&\\
q(t;{\bf x},{\bf y})\sim\frac{C_{2}({\bf x},{\bf 0})(\rho_{2}({\bf y})-J(0;{\bf y}))}{t\ln^{2}{t}},\quad &d=2,&\\
q(t;{\bf x},{\bf y})\sim\frac{C_{d}({\bf x},{\bf 0})(\rho_{d}({\bf
y})-J(0;{\bf y}))}{t^{d/2}},\quad &d\geq3,&
\end{eqnarray*}
where $C_{1}({\bf x},{\bf y})-C_{1}({\bf x},{\bf 0})J(0;{\bf y})>0$
for all ${\bf x},{\bf y}\in\mathbb{Z}$ and $J(0;{\bf
y})<\rho_{d}({\bf y})$ for $d\geq 2$ and ${\bf y}\in\mathbb{Z}^{d}$.
\end{Th}

To prove Theorem \ref{T-2} we apply the H\"{o}lder inequality
combined with results on connection between fractional moments of
random variables and fractional derivatives of their Laplace
transforms (see, e.g., \cite{Wolfe_75}). Theorem \ref{T-3} can be
considered as a corollary of Theorem \ref{T-2}.

\begin{Th}\label{T-3}
Let ${\sf E}{\xi}<1+h_{d}\alpha^{-1}(1-\alpha)$ and ${\sf E}{\xi^{1+\delta}}$
be finite for some $\delta\in(0,1]$. Then, for fixed
${\bf x},{\bf y}\in\mathbb{Z}^{d}$ and each $s\in[0,1]$, the following
equalities are valid
\begin{eqnarray*}
\lim\limits_{t\to\infty}{{\sf E}_{\bf x}\left(\left.s^{\mu(t;{\bf
y})}\right|\mu(t;{\bf y})>0\right)}= \frac{s\,C_{1}({\bf x},{\bf
y})-C_{1}({\bf x},{\bf 0})\left(J(0;{\bf y})-J(s;{\bf
y})\right)}{C_{1}({\bf x},{\bf y})-C_{1}({\bf x},{\bf 0})J(0;{\bf
y})},\quad &d=1,&\\
\lim\limits_{t\to\infty}{{\sf E}_{\bf x}\left(\left.s^{\mu(t;{\bf
y})}\right|\mu(t;{\bf y})>0\right)}=\ \frac{s\,\rho_{d}({\bf
y})-\left(J(0;{\bf y})-J(s;{\bf y})\right)}{\rho_{d}({\bf
y})-J(0;{\bf y})},\quad &d\geq2.&
\end{eqnarray*}
\end{Th}

{\sc Remark.} Comparing formulations of Theorems
\ref{T-1}--\ref{T-3} and the results on local particles numbers in
\emph{critical} CBRW on $\mathbb{Z}^{2}$ (see, e.g.,
\cite{TVP_Bulinskaya} and \cite{JOTP_Bulinskaya}) the following
conclusion suggests itself. Local particles numbers in critical CBRW
demonstrate ``subcritical'' behavior in the case of random walk on
integer lattice. Moreover, the scheme of proofs of Theorems
\ref{T-1}--\ref{T-3} is easily carried over to the case of critical
CBRW on $\mathbb{Z}^{2}$. Consequently, we may state that the
results of \cite{TVP_Bulinskaya} and \cite{JOTP_Bulinskaya},
concerning the local numbers in critical CBRW on $\mathbb{Z}^{2}$,
are valid under less restrictive conditions on the moments of
offspring number. Namely, it is sufficient to require finiteness of
${\sf E}{\xi}^{1+\delta}$ for some $\delta\in(0,1]$ instead of
condition ${\sf E}{\xi}^{2}<\infty$.

Concluding the first part of the paper let us note the close
relation between CBRW and superprocesses, namely, catalytic
super-Brownian motion with a single point of catalysis (see, e.g.,
\cite{Vatutin_Xiong_07} and references therein). It is of interest
that in view of \cite{VTY} CBRW may be considered as a queueing
system with a random number of independent servers. This gives an
opportunity for wide applications of the established results.

\section{Proof of Theorem \ref{T-1}}

Similarly to the proof of Theorem 1 in \cite{JOTP_Bulinskaya}, we
use backward and forward integral equations for the family of
functions $\{m(\cdot;{\bf x},{\bf y})\}_{{\bf x},{\bf
y}\in\mathbb{Z}^{d}}$. These equations coincide with equations (8)
and (9) in \cite{JOTP_Bulinskaya}, derived for the mean local
particles numbers in \emph{critical} CBRW on $\mathbb{Z}^{d}$, upon
the replacement of critical value
${\beta_{c}=(1-\alpha)a^{-1}G^{-1}_{0}({\bf 0},{\bf 0})}$ by value
$\beta$, that is,
\begin{eqnarray}
m(t;{\bf x},{\bf y})=p(t;{\bf x},{\bf
y})&+&\left(1-\frac{a}{1-\alpha}\right)\int\nolimits_{0}^{t}{p(t-u;{\bf
x},{\bf 0})m'(u;{\bf 0},{\bf
y})\,du}\nonumber\\
&+&\frac{a\beta}{1-\alpha}\int\nolimits_{0}^{t}{p(t-u;{\bf x},{\bf
0})m(u;{\bf 0},{\bf y})\,du},\label{m(t;x,y)backward}\\
m(t;{\bf x},{\bf y})=p(t;{\bf x},{\bf
y})&+&\left(\frac{1-\alpha}{a}-1\right)\int\nolimits_{0}^{t}{m(u;{\bf
x},{\bf 0})p'(t-u;{\bf 0},{\bf
y})\,du}\nonumber\\
&+&\beta\int\nolimits_{0}^{t}{m(u;{\bf x},{\bf 0})p(t-u;{\bf 0},{\bf
y})\,du}.\label{m(t;x,y)forward}
\end{eqnarray}
Relations \eqref{m(t;x,y)backward} and \eqref{m(t;x,y)forward}, as
well as equations (8) and (9) in \cite{JOTP_Bulinskaya}, are
obtained by means of the variation of constants formula applied to
backward and forward differential equations (5) and (6) in
\cite{JOTP_Bulinskaya} established in Banach space
$l_{\infty}(\mathbb{Z}^{d})$.

We also need the following auxiliary statement proved
similarly to Lemma 3.3.5 in \cite{Yarovaya_book} and Lemma 1
in \cite{JOTP_Bulinskaya}.

\begin{Lm}\label{L-monotonicity}
For each ${\bf y}\in\mathbb{Z}^{d}$ the function $m(t;{\bf y},{\bf y})$
does not increase in variable $t$.
\end{Lm}

Now we turn directly to the proof of Theorem \ref{T-1}. At first let us consider the case ${{\bf x}={\bf y}={\bf 0}}$. Apply
the Laplace transform to both sides of equality
\eqref{m(t;x,y)forward} and use the obtained relation to express function
$\widehat{m}(\lambda):=\int\nolimits_{0}^{\infty}{e^{-\lambda
t}m(t;{\bf 0},{\bf 0})\,dt}$, $\lambda\geq0$. Then
\begin{equation}\label{m(lambda)=}
\widehat{m}(\lambda)=\frac{G_{\lambda}({\bf 0},{\bf
0})}{1-\left((1-\alpha)a^{-1}-1\right)\int\nolimits_{0}^{\infty}{e^{-\lambda
t}p\,'(t;{\bf 0},{\bf 0})\,dt}-\beta G_{\lambda}({\bf 0},{\bf 0})}.
\end{equation}
Differentiate each side of the last relation in $\lambda$.
Consequently, taking into account the identity
${\int\nolimits_{0}^{\infty}{e^{-\lambda t}p\,'(t;{\bf 0},{\bf
0})\,dt}=\lambda G_{\lambda}({\bf 0},{\bf 0})-1}$ one has
\begin{equation}\label{m'(lambda)}
\widehat{m}'(\lambda)=\frac{(1-\alpha)a^{-1}G^{\,'}_{\lambda}({\bf
0},{\bf 0})+((1-\alpha)a^{-1}-1) G^{2}_{\lambda}({\bf 0},{\bf
0})}{\left((1-\alpha)a^{-1}-\left((1-\alpha)a^{-1}-1\right)\lambda
G_{\lambda}({\bf 0},{\bf 0})-\beta G_{\lambda}({\bf 0},{\bf
0})\right)^{2}}.
\end{equation}
According to Tauberian Theorem 2 in \S 5 of Chapter XIII in
\cite{Feller} along with Corollary 43 in \cite{V_Lectures}, relation
\eqref{p(t;x,y)sim} implies
\begin{eqnarray*}
G_{\lambda}({\bf 0},{\bf
0})\sim\frac{\gamma_{1}\sqrt{\pi}}{\sqrt{\lambda}},\quad
G^{\,'}_{\lambda}({\bf 0},{\bf
0})\sim-\frac{\gamma_{1}\sqrt{\pi}}{2\lambda^{3/2}},\quad &d=1,&\\
G_{\lambda}({\bf 0},{\bf
0})\sim\gamma_{2}\ln{\frac{1}{\lambda}},\quad G^{\,'}_{\lambda}({\bf
0},{\bf 0})\sim-\frac{\gamma_{2}}{\lambda},\quad &d=2,&
\end{eqnarray*}
as $\lambda\to0+$. Substituting these asymptotic equalities into
\eqref{m'(lambda)} for $d=1$ and $d=2$, respectively, we find
\begin{eqnarray*}
\widehat{m}'(\lambda)\sim-\frac{1-\alpha}{2\,a\,\gamma_{1}\,\sqrt{\pi}\,\beta^{2}\,\sqrt{\lambda}},\quad
&d=1,&\\
\widehat{m}'(\lambda)\sim-\frac{1-\alpha}{a\,\gamma_{2}\,\beta^{2}\,\lambda\,\ln^{2}{\lambda}},\quad
&d=2,&
\end{eqnarray*}
as $\lambda\to0+$. Applying Corollary 43 in \cite{V_Lectures} to the
above relations we come to the assertion of Theorem~\ref{T-1} when
${\bf x}={\bf y}={\bf 0}$ and $d=1$ or $d=2$.

For $d\geq3$, we employ another approach, namely, the representation
of complex-valued measures in terms of Banach algebras. In view of
\eqref{p(t;x,y)sim} two cumulative distribution functions, having
the Laplace transforms $G_{\lambda}({\bf 0},{\bf 0})$ and
$\int\nolimits_{0}^{\infty}{e^{-\lambda t}p'(t;{\bf 0},{\bf
0})\,dt}$, respectively, possess tails equivalent to constants (the second constant being
zero) multiplied by the same function $t^{1-d/2}$. Hence, due to
Lemma~6 in \cite{VT_Siberia} and formula \eqref{m(lambda)=} the
following asymptotic equality holds true
$$\int\nolimits_{t}^{\infty}{m(u;{\bf 0},{\bf 0})\,du}\sim\frac{2\,(1-\alpha)\,a\,\gamma_{d}}{(d-2)(1-\alpha-
a\,\beta\,G_{0}({\bf 0},{\bf 0}))^{2}\,t^{d/2-1}},\quad
t\to\infty.$$ Whence by Lemma~\ref{L-monotonicity} and the classical
results on differentiating of asymptotic formulae (see, e.g.,
\cite{Bruijn}, Chapter~7, Section~3) we get the assertion of
Theorem~\ref{T-1} when ${\bf x}={\bf y}={\bf 0}$ and $d\geq3$.

Let us consider the case ${\bf x}\neq{\bf 0}$ and ${\bf y}={\bf 0}$.
Integration by parts permits to rewrite the family of equations
\eqref{m(t;x,y)backward} as follows
\begin{eqnarray}
m(t;{\bf x},{\bf 0})=\frac{a}{1-\alpha}\,p(t;{\bf x},{\bf
0})&+&\left(1-\frac{a}{1-\alpha}\right)\int\nolimits_{0}^{t}{m(t-u;{\bf
0},{\bf 0})\,p\,'(u;{\bf x},{\bf
0})\,du}\nonumber\\
&+&\frac{a\,\beta}{1-\alpha}\int\nolimits_{0}^{t}{p(t-u;{\bf
x},{\bf 0})\,m(u;{\bf 0},{\bf 0})\,du},\label{m(t;x,0)1}\\
\frac{a}{1-\alpha}\,m(t;{\bf 0},{\bf
0})=\frac{a}{1-\alpha}\,p(t;{\bf 0},{\bf
0})&+&\left(1-\frac{a}{1-\alpha}\right)\int\nolimits_{0}^{t}{m(t-u;{\bf
0},{\bf 0})\,p\,'(u;{\bf 0},{\bf 0})\,du}\nonumber\\
&+&\frac{a\,\beta}{1-\alpha}\int\nolimits_{0}^{t}{p(t-u;{\bf 0},{\bf
0})\,m(u;{\bf 0},{\bf 0})\,du}.\label{m(t;x,0)2}
\end{eqnarray}
Subtracting equation \eqref{m(t;x,0)1} from \eqref{m(t;x,0)2} we come to
\begin{eqnarray}
\frac{a}{1-\alpha}\,m(t;{\bf 0},{\bf 0})-m(t;{\bf x},{\bf
0})&=&\frac{a}{1-\alpha}\left(p(t;{\bf 0},{\bf 0})- p(t;{\bf x},{\bf
0})\right)\nonumber\\
&+&\left(1-\frac{a}{1-\alpha}\right)\int\nolimits_{0}^{t}m(t-u;{\bf
0},{\bf 0})\left(p\,'(u;{\bf 0},{\bf 0})-p\,'(u;{\bf x},{\bf
0})\right)\,du\nonumber\\
&+&\frac{a\,\beta}{1-\alpha}\int\nolimits_{0}^{t}{m(t-u;{\bf 0},{\bf
0})\left(p(u;{\bf 0},{\bf 0})-p(u;{\bf x},{\bf
0})\right)\,du}.\label{m(t;x,0)-m(t;0,0)=}
\end{eqnarray}
Employing the results on differentiating of asymptotic formulae once
again (see, e.g., \cite{Bruijn}, Chapter~7, Section~3), as well as
relation \eqref{p(t;x,y)sim} and inequality $p\,''(t;{\bf 0},{\bf
0})\geq p\,''(t;{\bf x},{\bf 0})$, $t\geq0$, implied by (2.1.15) in
\cite{Yarovaya_book}, we find that
\begin{equation}\label{p'(t;0,0)-p'(t;x,0)sim}
p'(t;{\bf 0},{\bf 0})-p'(t;{\bf x},{\bf
0})\sim-\frac{(d+2)\,\widetilde{\gamma}_{d}({\bf
x})}{2\,t^{d/2+2}},\quad t\to\infty.
\end{equation}
Therefore, on account of Lemma~5.1.2 (``lemma on convolutions'') in
\cite{Yarovaya_book} along with formula \eqref{p(t;x,y)sim} and the
proved part of Theorem~\ref{T-1} we deduce from
\eqref{m(t;x,0)-m(t;0,0)=} that
\begin{eqnarray}
m(t;{\bf x},{\bf 0})&\sim&m(t;{\bf 0},{\bf
0})\left(1-\frac{a\,\beta}{1-\alpha}\int\nolimits_{0}^{\infty}
{\left(p(u;{\bf 0},{\bf 0})-p(u;{\bf x},{\bf
0})\right)\,du}\right)\nonumber\\
&-&\frac{a\,\widetilde{\gamma}_{d}({\bf
x})}{(1-\alpha)\,t^{d/2+1}}\left(1+\beta\int\nolimits_{0}^{\infty}{m(u;{\bf
0},{\bf 0})\,du}\right),\quad t\to\infty.\label{m(t;x,0)sim}
\end{eqnarray}
However by virtue of already established part of Theorem~\ref{T-1}
we see that $\int\nolimits_{0}^{\infty}{m(u;{\bf 0},{\bf
0})\,du}<\infty$ for all $d\in\mathbb{N}$ and, moreover, according
to \eqref{m(lambda)=} for $\lambda=0$ one gets
\begin{equation}\label{int_m(t;0,0)}
\int\nolimits_{0}^{\infty}{m(u;{\bf 0},{\bf
0})\,du}=-\beta^{-1}\quad\mbox{if}\quad d=1\quad\mbox{or}\quad d=2.
\end{equation}
Hence we conclude that only the first summand in the right-hand side
of \eqref{m(t;x,0)sim} makes contribution to the asymptotic behavior
of $m(t;{\bf x},{\bf 0})$. So, the assertion of Theorem~\ref{T-1}
for ${\bf x}\neq{\bf 0}$ and ${\bf y}={\bf 0}$ is entailed by
relation \eqref{m(t;x,0)sim} and the proved part of this theorem
when ${\bf x}={\bf y}={\bf 0}$.

Let now ${\bf x}\in\mathbb{Z}^{d}$ and ${\bf y}\neq{\bf 0}$. In view
of \eqref{m(t;x,y)forward} one has
\begin{eqnarray*}
m(t;{\bf x},{\bf 0})-m(t;{\bf x},{\bf y})&=&p(t;{\bf x},{\bf
0})-p(t;{\bf x},{\bf
y})\\
&+&\left(\frac{1-\alpha}{a}-1\right)\int\nolimits_{0}^{t}{m(t-u;{\bf
x},{\bf 0})\left(p\,'(u;{\bf 0},{\bf 0})-p\,'(u;{\bf 0},{\bf
y})\right)du}\\
&+&\beta\int\nolimits_{0}^{t}{m(t-u;{\bf x},{\bf 0})\left(p(u;{\bf
0},{\bf 0})-p(u;{\bf 0},{\bf y})\right)du}.
\end{eqnarray*}
Then taking into account formulae \eqref{p(t;x,y)sim} and
\eqref{p'(t;0,0)-p'(t;x,0)sim} along with Lemma 5.1.2 in
\cite{Yarovaya_book} and the proved part of Theorem \ref{T-1} we
find that
\begin{eqnarray}
m(t;{\bf x},{\bf y})&\sim& m(t;{\bf x},{\bf
0})\left(\frac{1-\alpha}{a}-
\beta\int\nolimits_{0}^{\infty}{\left(p(u;{\bf 0},{\bf 0})-p(u;{\bf
0},{\bf y})\right)du}\right)\nonumber\\
&+&\frac{\widetilde{\gamma}_{d}({\bf x})-\widetilde{\gamma}_{d}({\bf
y}-{\bf x})}{t^{d/2+1}}-\frac{\beta\,\widetilde{\gamma}_{d}({\bf
y})}{t^{d/2+1}}\int\nolimits_{0}^{\infty}{m(u;{\bf x},{\bf
0})\,du},\quad t\to\infty.\label{m(t;x,y)sim}
\end{eqnarray}
By the established part of Theorem \ref{T-1} applied one again we
come to inequality $\int\nolimits_{0}^{\infty}{m(u;{\bf x},{\bf
0})\,du}<\infty$. Moreover, similarly to the verification of equality
\eqref{int_m(t;0,0)} we check that
$\int\nolimits_{0}^{\infty}{m(u;{\bf x},{\bf 0})\,du}=-\beta^{-1}$
for $d=1$ or $d=2$. Thus, the statement of Theorem \ref{T-1} for
${\bf x}\in\mathbb{Z}^{d}$ and ${\bf y}\neq{\bf 0}$ is implied by
relation \eqref{m(t;x,y)sim} and the proved part of
Theorem~\ref{T-1} for ${\bf x}\in\mathbb{Z}^{d}$ and ${\bf y}={\bf
0}$.

To complete the proof of Theorem~\ref{T-1} one has to make sure only
that functions $C_{d}(\cdot,\cdot)$, $d\in\mathbb{N}$, are strictly
positive. It is easy except for the case $d=1$ when ${\bf x}\neq{\bf
0}$ and ${\bf y}\neq{\bf 0}$. Let us show that in this case
$C_{1}({\bf x},{\bf y})>0$ as well. For this purpose we turn to
function $H_{{\bf x},{\bf 0}}(t)$, $t\geq0$, which is a cumulative
distribution function of time from leaving point ${\bf x}$ till the
first hitting point ${\bf 0}$ in the framework of the random walk
generated by matrix $A$. Clearly,
\begin{equation}\label{C1(x,y)>0}
p(t;{\bf x},{\bf y})-\int\nolimits_{0}^{t}{p(t-u;{\bf 0},{\bf y})\,d
H_{{\bf x},{\bf 0}}(u)}\geq0,\quad t\geq0.
\end{equation}
Then by virtue of the evident identity $p(t;{\bf x},{\bf
0})=\int\nolimits_{0}^{t}{p(t-u;{\bf 0},{\bf 0})\,d H_{{\bf x},{\bf
0}}(u)}$ one has
\begin{eqnarray*}
p(t;{\bf x},{\bf y})-\int\nolimits_{0}^{t}{p(t-u;{\bf 0},{\bf y})\,d
H_{{\bf x},{\bf 0}}(u)}&=&p(t;{\bf x},{\bf y})-p(t;{\bf x},{\bf
0})\\
&+&\int\nolimits_{0}^{t}{\left(p(t-u;{\bf 0},{\bf 0})-p(t-u;{\bf
0},{\bf y})\right)\,d H_{{\bf x},{\bf 0}}(u)}.
\end{eqnarray*}
Since with the help of relation \eqref{p(t;x,y)sim}, Lemma 5.1.2 in
\cite{Yarovaya_book} and Lemma 3 in \cite{LMJ_Bulinskaya} one can
find the asymptotic behavior of the right-hand side of the latter
equality, we establish that, as $t\to\infty$,
\begin{eqnarray*}
p(t;{\bf x},{\bf y})-\int\nolimits_{0}^{t}{p(t-u;{\bf 0},{\bf y})\,d
H_{{\bf x},{\bf 0}}(u)}&\sim&\frac{\widetilde{\gamma}_{1}({\bf
x})+\widetilde{\gamma}_{1}({\bf
y})-\widetilde{\gamma}_{1}({\bf y}-{\bf x})}{t^{3/2}}\\
&+&\frac{(1-\alpha-a\,\rho_{1}({\bf x}))(1-\alpha-a\,\rho_{1}({\bf
y}))}{2\,a^{2}\,\gamma_{1}\,\pi\,\beta^{2}\,t^{3/2}}.
\end{eqnarray*}
Hence, it follows from \eqref{C1(x,y)>0} that
$$\frac{\rho_{1}({\bf x})\rho_{1}({\bf y})}{2\,\gamma_{1}\,\pi\,\beta^{2}}+\widetilde{\gamma}_{1}({\bf x})+
\widetilde{\gamma}_{1}({\bf y})-\widetilde{\gamma}_{1}({\bf y}-{\bf
x})\geq\frac{1-\alpha}{2\,a\,\gamma_{1}\,\pi\,\beta^{2}}\,\left(\rho_{1}({\bf
x})+\rho_{1}({\bf y})-\frac{1-\alpha}{a}\right).$$ However the
left-hand side of this inequality appears to be $C_{1}({\bf x},{\bf
y})$, whereas the right-hand side is strictly positive since
$\rho_{1}({\bf z})>(1-\alpha)a^{-1}$ for ${\bf z}\neq{\bf 0}$. In
its turn, the last inequality is satisfied due to definition of
function $\rho_{1}(\cdot)$ and negativeness of $\beta$ in
subcritical regime for $d=1$.

Therefore, Theorem~\ref{T-1} is proved completely. $\square$

\section{Proofs of Theorems~\ref{T-2} and \ref{T-3}}

It is not difficult to verify (following the scheme in
\cite{JOTP_Bulinskaya}) that in subcritical CBRW on
$\mathbb{Z}^{d}$, as well as in critical CBRW, for all ${\bf x},{\bf
y}\in\mathbb{Z}^{d}$, $s\in[0,1]$ and $t\geq0$, the non-linear
integral equations hold true
\begin{equation}\label{q(s,t;x,y)_equation}
q(s,t;{\bf x},{\bf y})=(1-s)m(t;{\bf x},{\bf
y})-\int\nolimits_{0}^{t}{m(t-u;{\bf x},{\bf 0})\Phi(q(s,u;{\bf
0},{\bf y}))\,du}.
\end{equation}
So, the following upper estimate for $q(s,t;{\bf x},{\bf y})$
ensues in view of non-negativeness of functions $m(\cdot;{\bf x},{\bf 0})$
and $\Phi(\cdot)$
\begin{equation}\label{q(s,t;x,y)_upper_bound}
q(s,t;{\bf x},{\bf y})\leq(1-s)m(t;{\bf x},{\bf y}).
\end{equation}

\begin{Lm}\label{L-fractional}
If ${\sf E}{\xi}<1+h_{d}\alpha^{-1}(1-\alpha)$ and ${\sf
E}{\xi^{1+\delta}}<\infty$ for $\delta\in(0,1]$ then for some
positive constants $K_{1}$ and $K_{2}$ the inequalities are valid
\begin{eqnarray}
\Phi(s)&\leq&K_{1}s^{1+\delta},\quad s\in[0,1],\label{Phi(s)<=}\\
{\sf E}_{\bf x}{\mu(t;{\bf y})^{1+\delta}}&\leq&K_{2}\,m(t;{\bf
x},{\bf y}),\quad t\geq t_{0}({\bf x},{\bf y}),\;{\bf x},{\bf
y}\in\mathbb{Z}^{d}\label{E_mu(t;y)_delta<=}
\end{eqnarray}
with a certain non-negative function $t_{0}(\cdot,\cdot)$.
\end{Lm}
{\sc Proof.} At first we consider the case $0<\delta<1$. Let us take
advantage of the connection between fractional moments of random
variables and fractional derivatives of their Laplace transforms.
For the first time such results were obtained in \cite{Wolfe_75}
where the traditional notion of Riemann-Liouville fractional
derivative was used. However it is more convenient to involve the
up-to-date counterpart of these results, namely, Lemma~2.1 in
\cite{Klar_2003}, which gives
\begin{equation}\label{E_xi_1+delta=}
{\sf
E}{\xi^{1+\delta}}=\frac{\delta(1+\delta)}{\Gamma(1-\delta)}\int\nolimits_{0}^{\infty}
{\frac{\Phi(1-e^{-v})+\alpha\,f'(1)(e^{-v}-1+v)}{\alpha\,v^{2+\delta}}\,dv}.
\end{equation}
In view of finiteness of ${\sf E}{\xi}^{1+\delta}$ the latter
equality entails
$\int\nolimits_{0}^{1}{v^{-2-\delta}\Phi(v)\,dv}<\infty.$ Taking into account that
$\Phi'(s)\geq0$, $s\in[0,1]$, integration by parts implies relation \eqref{Phi(s)<=} for
$\delta\in(0,1)$.

Turn to verification of \eqref{E_mu(t;y)_delta<=} when $0<\delta<1$.
Applying Lemma~2.1 in \cite{Klar_2003} once again we come to the
following equality
$${\sf E}_{\bf x}{\mu(t;{\bf y})^{1+\delta}}=\frac{\delta(1+\delta)}{\Gamma(1-\delta)}
\int\nolimits_{0}^{\infty}{\frac{v\,m(t;{\bf x},{\bf
y})-q(e^{-v},t;{\bf x},{\bf y})}{v^{2+\delta}}\,dv}.$$ Substituting
formula \eqref{q(s,t;x,y)_equation} into the latter relation we get
\begin{eqnarray}
{\sf E}_{\bf x}{\mu(t;{\bf y})^{1+\delta}}&=&m(t;{\bf x},{\bf
y})\frac{\delta(1+\delta)}{\Gamma(1-\delta)}
\int\nolimits_{0}^{\infty}{\frac{e^{-v}-1+v}{v^{2+\delta}}\,dv}\nonumber\\
&+&\frac{\delta(1+\delta)}{\Gamma(1-\delta)}\int\nolimits_{0}^{\infty}
{\frac{1}{v^{2+\delta}}\int\nolimits_{0}^{t}{m(t-u;{\bf x},{\bf
0})\Phi(q(e^{-v},u;{\bf 0},{\bf
y}))\,du}\,dv}.\label{fractional_moment_1}
\end{eqnarray}
Obviously, in equality \eqref{fractional_moment_1} the first
integral converges. Let us estimate the double integral. According
to \eqref{q(s,t;x,y)_upper_bound} and inequality $\Phi(\kappa
s)\leq\kappa \Phi(s)$, $\kappa,s\in[0,1]$, guaranteed by
convexity property of function $\Phi(\cdot)$, one has
$$\int\limits_{0}^{\infty}
{\frac{1}{v^{2+\delta}}\int\limits_{0}^{t}{m(t-u;{\bf x},{\bf
0})\Phi(q(e^{-v},u;{\bf 0},{\bf
y}))\,du}\,dv}\leq\int\limits_{0}^{t}{m(t-u;{\bf x},{\bf 0})m(u;{\bf
0},{\bf y})\,du}\int\limits_{0}^{\infty}
{\frac{\Phi(1-e^{-v})}{v^{2+\delta}}\,dv}.$$ At the right-hand side
of the above inequality the integral in variable $u$ is equivalent (up to a constant
factor) to the function $m(t;{\bf x},{\bf y})$, as $t\to\infty$, on account of Theorem~\ref{T-1} and
Lemma~5.1.2 in \cite{Yarovaya_book}. Furthermore, in the same
inequality the integral in variable $v$ converges by virtue of
formula \eqref{E_xi_1+delta=} combined with finiteness of ${\sf
E}{\xi}^{1+\delta}$. Therefore, this
arguing together with relation \eqref{fractional_moment_1} lead to
\eqref{E_mu(t;y)_delta<=} with $0<\delta<1$.

Now we consider the case $\delta=1$. The identity $f''(1)={\sf
E}{\xi(\xi-1)}$ and the Lemma conditions imply the existence of
$f''(1)$ and, consequently, the existence of $\Phi''(0)$. Moreover,
since $\Phi(0)=0$ and $\Phi'(0)=0$, inequality \eqref{Phi(s)<=} for
$\delta=1$ is proved. Let us take the second left derivatives at $s=1$
for each side of equality \eqref{q(s,t;x,y)_equation}.
Using relations $m(t;{\bf x},{\bf y})=-\left.\partial_{s} q(s,t;{\bf
x},{\bf y})\right|_{s=1}$ and ${\sf E}_{\bf x}{\mu(t;{\bf
y})(\mu(t;{\bf y})-1)}=-\left.\partial^{\,2}_{s s}\,q(s,t;{\bf
x},{\bf y})\right|_{s=1}$, we obtain
$${\sf
E}_{\bf x}{\mu(t;{\bf y})(\mu(t;{\bf y})-1)}=\alpha
f''(1)\int\nolimits_{0}^{t}{m(t-u;{\bf x},{\bf 0})\left(m(u;{\bf
0},{\bf y})\right)^{2}\,du}.$$ In accordance with Theorem~\ref{T-1}
and  Lemma~5.1.2 in \cite{Yarovaya_book}, in the latter equality the
integral behaves as $m(t;{\bf x},{\bf 0})$ up to a constant factor,
as $t\to\infty$. Hence, on account of Theorem~\ref{T-1} this entails
the desired inequality \eqref{E_mu(t;y)_delta<=} when $\delta=1$.
Lemma~\ref{L-fractional} is proved completely. $\square$

Let us turn to proving Theorem~\ref{T-2}. It is easily
verified that due to formulae \eqref{q(s,t;x,y)_upper_bound} and
\eqref{Phi(s)<=} along with Theorem~\ref{T-1} and the proof scheme
of Lemma~4 in \cite{TVP_Bulinskaya}, one gets
\begin{equation}\label{int_convolution_sim}
\int\nolimits_{0}^{t}{m(t-u;{\bf x},{\bf 0})\Phi(q(s,u;{\bf 0},{\bf
y}))\,du}\sim m(t;{\bf x},{\bf 0}) J(s;{\bf y}),\quad t\to\infty.
\end{equation}
Observe that the function $J(\cdot;\cdot)$, appearing in formulations of Theorems~\ref{T-2} and
\ref{T-3}, is defined
correctly in view of upper estimate \eqref{q(s,t;x,y)_upper_bound}
for $q(s,t;{\bf x},{\bf y})$. Let us find the lower estimate for
$q(t;{\bf x},{\bf y})$. By the H\"{o}lder inequality one has
$${\sf E}_{\bf x}{\mu(t;{\bf y})}={\sf E}_{\bf x}{\mu(t;{\bf y})\mathbb{I}(\mu(t;{\bf y})>0)}\leq
\left({\sf E}_{\bf x}{\mu(t;{\bf
y})^{1+\delta}}\right)^{1/(1+\delta)}\left({\sf E}_{\bf
x}{\mathbb{I}(\mu(t;{\bf
y})>0)^{(1+\delta)/\delta}}\right)^{\delta/(1+\delta)}$$ where
$\mathbb{I}(B)$ denotes the indicator of set $B$. Rewrite the last
inequality as follows
$$q(t;{\bf x},{\bf y})\geq\frac{\left(m(t;{\bf x},{\bf y})\right)^{(1+\delta)/\delta}}
{\left({\sf E}_{\bf x}{\mu(t;{\bf
y})^{1+\delta}}\right)^{1/\delta}},\quad t>0.$$ Now employing
\eqref{E_mu(t;y)_delta<=} we come to the desired lower estimate for
$q(t;{\bf x},{\bf y})$
\begin{equation}\label{q(t;x,y)_lower_bound}
q(t;{\bf x},{\bf y})\geq K_{2}^{\,-1/\delta}m(t;{\bf x},{\bf
y}),\quad t\geq t_{0}({\bf x},{\bf y}),\;{\bf x},{\bf
y}\in\mathbb{Z}^{d}.
\end{equation}
Combining formulae \eqref{q(s,t;x,y)_equation} and
\eqref{int_convolution_sim}, when $s=0$, with estimate
\eqref{q(t;x,y)_lower_bound} we conclude that, as $t\to\infty$,
$$q(t;{\bf x},{\bf y})\sim m(t;{\bf x},{\bf y})-m(t;{\bf x},{\bf 0})J(0;{\bf y})\quad\mbox{and}\quad
J(0;{\bf y})<\lim\limits_{t\to\infty}{\frac{m(t;{\bf x},{\bf
y})}{m(t;{\bf x},{\bf 0})}}$$ for fixed ${\bf x},{\bf
y}\in\mathbb{Z}^{d}$. Just these relations amount to validity of
Theorem~\ref{T-2}. $\square$

Let us prove Theorem~\ref{T-3}. Applying formulae
\eqref{q(s,t;x,y)_equation} and \eqref{int_convolution_sim} once
again we have
$$q(s,t;{\bf x},{\bf y})\sim(1-s)m(t;{\bf x},{\bf y})-m(t;{\bf x},{\bf 0})J(s;{\bf y}),\quad t\to\infty,
\;{\bf x},{\bf y}\in\mathbb{Z}^{d}.$$ Then with the help of the identity
$$\lim\limits_{t\to\infty}{\sf E}_{\bf x}{\left(\left.s^{\mu(t;{\bf y})}\right|\mu(t;{\bf y})>0\right)}=
1-\lim\limits_{t\to\infty}{\frac{q(s,t;{\bf x},{\bf y})}{q(t;{\bf
x},{\bf y})}}$$ and Theorem~\ref{T-2} we obtain the assertion of
Theorem~\ref{T-3}. $\square$

\end{document}